\documentclass[11pt,pdftex]{amsart}
\usepackage{graphicx}
\usepackage{amssymb}
\usepackage{amsmath}

\newcommand{\pz}{P(z,\bar z)}
\newcommand{\g}{\gamma}

\newcommand{\ab}[1]{\vert z\vert^{#1}}
\newcommand{\cdva}{{\mathbb C^2}}

\newcommand{\te}{\theta}

\newcommand{\zz}{(z,\bar z)}

\newcommand{\fz}{{F(z,\bar z,u)}}
\newcommand{\fzs}{{\textstyle{F^*(z^*,\bar z^*,u^*)}}}

\newcommand{ \al}{\alpha}
\newcommand{\ppm}{m'}

\newcommand{\tf}{\tilde F}

\begin{document}
\title{
 Normal forms for hypersurfaces of finite type in
 $ \mathbb C^2$}
\author{Martin Kol\'a\v r}
\address{Department of Mathematical Analysis, Masaryk University,
\newline Janackovo nam. 2a, 662 95 Brno } \email {mkolar@math.muni.cz }
\maketitle

\begin{abstract}
We  construct normal forms  for Levi degenerate hypersurfaces of
finite type in $\mathbb C^2$. As one consequence, an explicit
solution   to the problem of local biholomorphic  equivalence   is
obtained. Another consequence
  determines the dimension of the stability group of the hypersurface.
\end{abstract}

\section{Introduction}

Levi degenerate hypersurfaces have been intensively studied  since the pioneering work
of J. J. Kohn ([K]), which introduced the concept of finite type.
 On the one hand, there is now deep qualitative theory which links
local geometry and analysis
on pseudoconvex domains of finite type (see [CS] for references).
On the other hand,
some fundamental quantitative results from the
nondegenerate case still did not find analogy (cf. [BFG]).

Our aim
 is to show that the construction of
 normal forms, developed by J. Moser
in  [CM],
has a natural generalization to Levi degenerate  hypersurfaces of
 finite type in  dimension two.
As a  consequence, the problem of local biholomorphic equivalence of
two hypersurfaces is
reduced  to a low dimensional algebraic calculation.
In fact, the remaining algebra behind
normal forms
is
much simpler than in the nondegenerate case,
as there are less  symmetries of the model hypersurface.
Another application gives precise information on the dimension of the stability
group of the hyperfurface at the given point.

Being the main motivation,
we first formulate the local equivalence problem.
Let $M_1, M_2 \subseteq \mathbb C^n$ be  real analytic hypersurfaces of real
codimension one and $p_1, p_2$ be points on $M_1$  and  $M_2$, respectively.
If
 there is an invertible holomorphic map $f$ defined in a
neighbourhood $U$ of
$p_1$
such that  $f(p_1) =p_2$ and $f(M_1 \cap U) \subseteq M_2$, the two hypersurfaces are said to be locally biholomorphically
equivalent.
The problem is to find a complete, computable set of local invariants,
which
provide  a way to recognize equivalent hypersurfaces.

The history of the problem starts with
  H. Poincar\'e, who showed that  the Riemann mapping theorem
has no equivalent in higher dimensions, and nontrivial local   invariants do  exist.
 For $n=2$, a solution was found
 by E. Cartan in [C1],[C2],  under the assumption that the hypersurfaces
are regular, i.e., the Levi form is nondegenerate at $p_1, p_2$.
For $n\ge 2$, with the same regularity assumption, the problem
was solved by Chern and Moser. While the result
 of
 the second part of [CM] generalizes
Cartan's result to higher dimensions (a result obtained independently
  by N. Tanaka),  the first part
gives
a different
solution, in terms of normal coordinates
and corresponding normal forms.

Both approaches in [CM] start with analysis of homogeneous model
hyperquadrics and their local symmetry groups.  In the
construction of normal forms the defining equation of the
hypersurface is put into a form in which certain  terms in its
Taylor expansion vanish. In effect, the hypersurface is osculated
by the model hyperquadric to a high order.  Choosing the
appropriate vanishing condition one can achieve that the transformation into
normal form is determined uniquely, up to the action of the local
symmetry group of the model hyperquadric. In the first step normal
forms are obtained by algebraic manipulation of formal power
series. The second step proves convergence of the transformation.

 The results of
Chern and Moser inspired a lot of subsequent work by many authors
(e.g. [F], [J], [V], [B], [W], [EHZ], [ELZ] and many others). Partial
results in constructing normal forms on Levi degenerate
 hypersurfaces   in $\cdva$
 were obtained by  Barletta-Bedford (hypersurfaces with a special symmetry), P. Wong
 (a subclass of hypersurfaces of type four) and  N. Stanton (rigid
 hypersurfaces).
P. Ebenfelt in [E]  constructed normal forms for a class of
hypersurfaces in $\mathbb C^n$ of finite type three.

The main difference, which seems  inherent to degenerate
hypersurfaces, is that normal forms are given by formal power
series which need not converge (although we do not give any
explicit example of a divergent normal form). The fact that a
normal form construction solves the local equivalence problem
relies on the
 essential result of M.S.Baouendi, P.Ebenfelt and L.P.Rothschild
 ([BER]),
 that any formal equivalence of two finite type hypersurfaces
has to converge.

 We  give the first step in the
construction of normal forms in Section 2. It is essentially a
partial normalization which removes low order harmonic terms. This
procedure reveals the type of the boundary point, an integer $k$,
and a relatively simple substitute for the model hyperquadric,
which plays a central role in Chern-Moser's theory. The
fundamental information contained in the model hypersurface is its
essential type, denoted by $l$. There are three qualitatively
different types of model hypersurfaces, which have to   be treated
separately.

In Section 3 we construct normal forms for the generic case, when
$l<\frac{k}2$ and the model is not a tube. First we consider
linear transformations and determine the local symmetry group of
the model hypersurface, whose real dimension is equal to one.
 Normal form
conditions are then defined (Definition 3.2) and  shown to
determine uniquely all coefficients of the biholomorphic
transformation (Proposition 3.3). The remaining cases are
considered in
 Section 4. Definition 4.1 describes normal forms for the case of the most
symmetric model hypersurface, when $l=\frac{k}2$. The symmetry
group is three dimensional in this case. Definition 4.3 gives
normal forms for the case when the model hypersurface is a tube.
Proposition 4.2 modifies Proposition 3.3 for these two cases.

Section 5 contains applications.
 First we show that the local equivalence problem is
solved by the construction of normal forms combined with the
convergence result of [BER] (Theorem 5.1). Then we apply our
results to obtain precise information about the dimension of the
stability group (Corollary 5.4).

  More applications of 
the results presented here are given in [Ko].

I would like to thank  M. Salah  Baouendi, Peter Ebenfelt and
Linda Preiss Rothschild for helpful comments on an earlier version
of this paper.

\section{Hypersurfaces of finite type}

Let $M \subseteq \mathbb C^2$ be a  real analytic hypersurface,
 and  $p$ be a point
on $ M$.
Let $r \in C^{\omega}$  be a local defining function, i.e.,
for  a neighbourhood $U$ of $p$

\begin{equation}
M \cap U = \{z \in U \ | \  r(z) = 0\},\tag{2.1}\end{equation}

and $\nabla r \neq 0$ in $ M \cap U$.

We recall a definition of a point of finite type.
For a smooth real valued function $f$ defined in neighbourhood
of $0$ in $\mathbb C$ let $\nu(f)$ denote the order of vanishing of $f$ at $0$.

{\bf Definition 2.1} $p$ is a point of finite type, if there
exists an integer $m$ such that
$$\nu (r\circ \gamma)\leq m$$
for all holomorphic maps $\g $  from a neighbourhood of $0\in \mathbb C$
into $\cdva$, satisfying $\g(0)=p$ and $\g^{\prime}(0)\neq 0$.
The smallest such integer is called the type of $p$.

Note that $M$ is Levi nondegenerate at $p$ if and only if $p$ is a
point of finite type two. This is the case considered by Cartan
and Chern-Moser. Hence, without any loss of generality,  we will
assume in the rest of this paper that $p$ is a point of finite
type $k$, where $k > 2$. Our aim is to assign to the pair $(M,p)$
formal power series in normal form. It will be unique up to the
action of the symmetry group of the model hypersurface, which will
be determined below.

Let $(z,w)$, $z = x + iy, w = u + iv$,  be local holomorphic
coordinates centered at $p$ such that the hyperplane $\{ v=0 \}$
is tangent to $M$ at $p$. Near $p$, by the implicit function
theorem,  $M$ is described as a graph of a function
\begin{equation} v= F(x,y,u), \tag{2.2} \end{equation}
where $F$ is a real valued function defined in a neighbourhood of
the origin in $\mathbb R^3$. Since $r$ is real analytic, $F$ is the
sum of its  Taylor expansion starting with 2-nd order terms, which
we express in terms of $(z, \bar z, u)$:
\begin{equation}\fz = \sum_{i+j+m\ge 2} a_{ijm} z^i \bar z^j u^m,\tag{2.3}
\end{equation}
where $a_{ijm}=\overline{a_{jim}}$.

We will consider
  holomorphic transformations

\begin{equation}
  z^*=z+ f(z,w), \ \ \ \ \   w^*=w+g(z,w),  \tag{2.4} \end{equation}
where $f$ and $g$ are represented by power series
\begin{equation} f(z,w)=\sum_{i,j =0}^{\infty} f_{ij} z^i w^j, \ \ \ \ \   g(z,w)= \sum_{i,j=0}^{\infty} g_{ij}
 z^i w^j.
 \tag{2.5} \end{equation}

Since we will  have to consider also formal hypersurfaces and
formal transformations, from now on
 we  allow both $F$ and $f,g$ to be formal power series.
In this case (2.2) - (2.5)  are interpreted in this  sense.

We are interested only in  transformations which  preserve the
above form, given by (2.2), (2.3).
 This will hold if and only if $f$ and $g$ contain no constant
term, and $\frac{\partial v^*}{\partial x}, \frac{\partial
v^*}{\partial y}, \frac{\partial v^*}{\partial u}$ are all zero at
the origin. In other words, we require that
\begin{equation}
f=0,\ \ g=0, \ \ g_z=0, \ \ Im\ g_w=0\ \ \  \text{at}\ \ \
z=w=0.\tag{2.6}
\end{equation}
 In the following, we will consider only
transformations satisfying (2.6). Let $F^*$ be the power series
describing  $M$ in new coordinates. Substituting (2.4) into
$v^*=\fzs$, we get the change of variables formula
\begin{equation}
  F^*(z+f, \bar z+\bar f,u+Re\ g)  = \fz + Im\ g(z,u+i\fz),
\tag{2.7} 
\end{equation}
 where $f$ and $Re\ g $ are also evaluated at
$(z,u+i\fz)$. It can be viewed as an equality of two power series
in $z, \bar z , u $ which allows to obtain  relations between the
 coefficients of $F^*$ and $F,f,g$.

The first step in our construction is the following standard
result.

{\bf Lemma 2.2. }
{\it There exist uniquely determined complex numbers
\newline $\al_2, \dots, \al_{k}, $ such that after the change of variable
\begin{equation}
w^* = w + \sum^{k}_{i=2} \al_i z^i,
\ \tag{2.8}
\end{equation}
 the defining equation (2.2) has form
\begin{equation}
v^* = \pz + o(\ab k, u^* ),\ \tag{2.9}
\end{equation}
where $P$ is a nonzero real valued homogeneous polynomial
of degree $k$
\begin{equation}
\pz = \sum_{j=1}^{k-1}
a_j z^j\bar z^{k-j}, \ \tag{2.10}
\end{equation}
 where $a_j \in \mathbb C$ and
$a_j = \overline{a_{k-j}}.$
}

{\it{proof. } }It follows from (2.7) that the value of $\al_j$
does not affect terms of order less than $j$ in $F^*$. We start
with second order terms and write $F$ as
$$ \fz = A_2z^2 + \bar A_2 \bar z^2 + a\vert z \vert^2 + o(\ab 2, u).$$
After a change of variable $w^* = w +  \sum^k_{i=2} \al_i z^i, \  $ we obtain
$$ F^*(z,\bar z, u^*) = A_2z^2 + \bar A_2 \bar z^2 + a\vert z \vert^2 +
 Im\;  \al_2 z^2
+o(\vert z \vert^2, u^*).$$ Here $a$ is the value of the Levi form
at $p$, so $a=0$. There is a unique  $\al_2$ which makes the
second order  terms vanish, namely $\al_2 = 2iA_2$. Now we proceed
by induction. Let for some $j >  2$ the coefficients  $\al_2,
\dots, \al_{j-1}$ be already determined, so that
$$ F^*(z,\bar z, u^*) = P_j(z, \bar z) + o(\ab j, u^*),$$
where $P_j(z, \bar z)$ is a real valued homogeneous polynomial of degree $j$.
Using Definition 2.1 it is easily verified  that if $j <k$,
 we must have  $P_j(z, \bar z) = Re\; A_j z^j$ for some
$A_j \in \mathbb C$. Hence we must take $\al_j = 2iA_j$. For $j=k$,
$\al_k$ is uniquely determined by the requirement that $P$ in
(2.10) contains no harmonic term. \

In order to preserve the form achieved by Lemma 2.1,  all
transformations  which we will consider have to satisfy
\begin{equation}
\frac{\partial^j g}{\partial z^j}=0 \ \  \text{at}\ \ \
z=w=0,\tag{2.11}
\end{equation}
 for $j=2,\dots k$, in addition to (2.6).

 Let $l$ denote the lowest
index in (2.10) for which $a_l \neq 0$. We have $1\leq l \leq
\frac{k}2$. Note that $l$ is the essential type of the model
hypersurface to $M$ at $p$, defined below.

The problem  now splits into three cases, depending on the form of
$P$.
 Two cases are "exceptional",  the case with  extra symmetries,
 when $2l = k$ and
$P = a_l\ab k$, and the  case when $P$ is equivalent to $(Re\;
z)^k$, which corresponds to a tube domain.  All other
hypersurfaces will be treated together, as the generic case. We
consider it first.

\section{Normal forms for generic  models}

As a next step we consider the effect of a linear transformation
\begin{equation}
 w^* = \delta w , \ \ \ \ \ \ \ z^* = \beta^{-1} z, \ \tag{3.1}
\end{equation}
where $\delta \in \mathbb R$ and $ \beta \in \mathbb C $. Part of (3.1)
will be used to normalize $P$, the other part will give
 the symmetry group of the model hypersurface.

In (2.10) we have for $j <\frac k2$
\begin{equation} 
 a_j z^j\bar z^{k-j} + a_{k-j} z^{k-j}\bar z^{j} = 2\ab {2j} Re\  a_j
z^{k-2j}.\ \tag{3.2}
\end{equation}
 We will denote $j' = k-2j$ for $1\leq
j\leq \frac k2$, considering $'$ as an operator which can be
applied to any integer from $1$ to $\frac{k}2$. In order to
normalize $P$ in the simplest possible way, we introduce the
following notation. Let $l=m_0 <m_1 <\dots <m_p <\frac k2$ be
 the indices in (2.10) for which $a_{m_i}\neq 0$. Denote by $L$ the greatest common divisor of
$\ppm_{0}, \dots, \ppm_{p}$ and let
$$q_i = \frac {\gcd(\ppm_0 , \dots, \ppm_i)}
{\gcd (\ppm_0, \dots,\ppm_{i+1})}$$
for $0 \leq i \leq p-1$.

{\bf {Lemma 3.1.}} 
{\it There exists $\beta \in \mathbb C$ such that after
the change of variables
\begin{equation}
 w^* = w , \ \ \ \ \ \ \ z^* = \beta^{-1} z, \ \tag{3.3}
\end{equation}
and dropping stars, $P$ satisfies
\begin{equation} 
a_l = 1 \ \tag{3.4}
\end{equation}
and
\begin{equation}
 \arg a_{m_{i+1}} \in [0,\frac{2\pi}{q_i})\ \tag{3.5}
\end{equation}
for $0 \leq i \leq p-1$.
These conditions  determine $P$ uniquely, while $\beta$ is unique
up to  multiplication by an $L$-th root of unity.
}

{\it{proof. }} The effect of (3.3) on (3.2) is simply
multiplication of the coefficient by $\vert \beta
\vert^{2j}\beta^{k-2j}$. The condition $a_l = 1$ determines
$\beta$ uniquely up to multiplication by an $l'$-th root of unity.
Each of the conditions (3.5) further reduces the number of
possible values of $\beta $ by a factor of $q_i$. Since  we have
$$ L \prod_{i=0}^{p-1} q_i = l',$$
there are $L$ possible values of $\beta$. On the other hand, if
$\beta $ is an $L$-th root of unity, then (3.3) preserves $P$, and
the conclusion follows. \

We will denote by $M_D$ the model hypersurface at $p$ :
$$M_D = \{(z,w) \in \cdva\ | \  v = \pz \},$$
which will be briefly called the model.

Up to now we have transformed  $M$ into form (2.9), (2.10), (3.4),
(3.5).
\newcommand  {\tplus}{+}
The following transformations preserve the model:
$$z^* =   \delta e^{i\theta} z,\ \  \ \ \   w^* = \delta^k w,$$
                           where
$e^{i\te}$ is an
$L$-th root  of unity and $\delta > 0$ for $k$ even or $\delta \in \mathbb R \setminus \{ 0 \}$ for $k$ odd.
Let $H$ denote the group of such transformations.
Hence $H= \mathbb R^+ \oplus {\mathbb Z}_L$ for $k$ even and $H= \mathbb R^* \oplus {\mathbb Z}_L$ for $k$ odd.

In the following we will assign weight $1$ to $z,\bar z$ and
weight $k$ to the variable $u$ in $F$ and $w$ in $f$ and $g$.
Hence a monomial $z^i \bar z^j u^m$ has weight $i+j+km$, and
$z^iw^j$ has weight $i+kj$. Using weights we can write (2.9) as
$$F = P + \text{terms of weight} \ge  k+1, $$
where $P$ satisfies (3.4) and (3.5).

We denote by $\mathcal F$ the set of formal power series of the
form
\begin{equation}
 \fz = \pz + \tf(z, \bar z, u),\tag{3.6} 
\end{equation}
where
$$\tf(z, \bar z, u) = \sum_{wt. > k} a_{ijm}z^i\bar z^j u^m.
$$
We  decompose $\tf$ into parts containing terms of equal weight:
$$\tf = \sum_{\nu=k+1}^{\infty} F_{\nu}.$$
We will also use  partial  expansion of $\tf$ in $z, \bar z$. Let
$$Z_{ij}(u) = \sum_{m} a_{ijm} u^m,$$
so we can write
$$\tf(z, \bar z, u) = \sum_{i,j}Z_{ij}(u)z^i\bar z^j.$$

Now we consider the group $\mathcal T_0$ of formal transformations
preserving this form. It is easily verified that $\mathcal T_0$
consists of transformations of the form
$$z^* =\delta e^{i\theta} z + \text{terms of weight} \ge 2,\ \  \ \ \   w^* = \delta^k w + \text{terms of weight}
\ge k+1,$$ where again $e^{i\te}$ is an $L$-th root  of unity and
$\delta > 0$ for $k$ even or $\delta \in \mathbb R \setminus \{ 0 \}$
for $k$ odd.

Let $\mathcal T$ be the set of formal transformations of the form
\begin{equation}
\begin{aligned}
  z^* = z + \sum_{wt. >1} & f_{ij} z^i w^j
\\
  w^* = w + \sum_{wt.>k} & g_{ij} z^i w^j.
\end{aligned} \tag{3.7}
\end{equation}

 Clearly, $\mathcal T$ is a group
under composition. Again we decompose the formal power series
 into parts of the same weight
$$f=\sum_{\nu = 2}^{\infty} f_{\nu} \ \ \ \ \
\text{and} \ \ \ \ \ g=\sum_{\nu =k+1}^{\infty} g_{\nu},$$ and
denote such an element of $\mathcal T$ by $(f,g)$.

If $F \in \mathcal F$ and $(f,g) \in \mathcal T$ it is easily
verified using (2.7) that the formal power series resulting from
transforming $F$ by $(f,g)$ is again in  $\mathcal F$. Hence
$\mathcal T$ acts on $\mathcal F$ via formula (2.7).

Now we check that any $\tau \in \mathcal T_0$ can be factored in a
unique way as
$$\tau = \phi \circ T,$$
with $\phi \in H$ and $ T \in \mathcal T$. Here $\phi$ is simply
the linear part of $\tau$.  Hence we can use elements of $H$ to
normalize transformations in $\mathcal T_0$ to satisfy
\begin{equation}
 f_z=0,\ \ \  Re\; g_w = 0 \ \ \ at\  z=w=0,\tag{3.8}
\end{equation}
i.e., to  be in $\mathcal T$.

In summary,  $\mathcal T$ consists precisely of transformations
satisfying normalization conditions (2.6), (2.11) and (3.8).

 For terms of weight $\mu>k$ in (2.7) we get from (2.7) and (3.6)
\begin{equation} 
\begin{aligned}  F^*_{\mu}(z, \bar z, u) +
2 Re \ P_z(z,\bar z)  f_{\mu-k+1}(z,u +&i\pz)= \\
 = F_{\mu}(z, \bar z,&u) + Im\ g_{\mu}(z,u+i\pz) + \dots
\end{aligned}
 \tag{3.9}
\end{equation}
where dots denote terms depending on $f_{\nu-k+1}, g_{\nu},
F_{\nu}, F^*_{\nu}$ for  $\nu < \mu$ , and $P_z = \frac{\partial
P}{\partial z}.$

The action of $\mathcal T$ defines an equivalence relation
 on $\mathcal F$ and our aim is to
find a condition which selects a unique element in each class of
equivalence. We will use the following
 scalar product on the vector space of homogeneous polynomials of degree $k-1$
without a harmonic term. If $Q = \sum^{k-2}_{j=1} \alpha_j z^j
\bar z^{k-1-j}$ and $S = \sum^{k-2}_{j=1} \beta_j z^j \bar
z^{k-1-j}$, then
$$(Q,S) = \sum^{k-2}_{j=1} \alpha_j \bar \beta_j.$$
This notation will be  used also for polynomials which may contain
a harmonic term, which is then ignored. We need this notation also
for polynomials whose coefficients depend on $u$. In particular,
we denote
\begin{equation}
(Z_{k-1}, P_z) = \sum_{j=1}^{k-2}Z_{j,k-1-j}  (j+1)\bar a_{j+1}.\tag{3.10}
\end{equation}

{\bf Definition 3.2.}
 We say that   $F$ is in  normal form if
\begin{equation} 
\begin{aligned} 
Z_{j0} & = 0, \ \ \ \ \ j=1,2,\dots,  \\
Z_{k-l+j,l} & = 0, \ \ \ \ \ j= 0,1,\dots, \\
Z_{2k-2l, 2l} & = 0, \\
(Z_{k-1}, P_z) & = 0.
\end{aligned}
  \tag{3.11} 
\end{equation}

We will prove

{\bf Proposition 3.3. }
{\it For any $F \in \mathcal F$ there is a
uniquely determined formal transformation $T \in \mathcal T$ which
transforms $F$ into normal form.
}

{\it{proof. }} By induction on weight we show that the condition
that $F^*$ satisfy (3.11) determines uniquely all coefficients of
$f$ and $g$ in (3.7). Let us consider terms of weight $\mu$ in
(2.7). For an analytic function $\phi$ of two variables we use the
identity
$$\phi(z, u+i P\zz) = \sum_{n=0}^{\infty} i^ {n}\frac{ \phi^{(n)}(z,u)}{n!}
P\zz^{n},$$
where $\phi^{(n)}$ denotes the $n$-th derivative of $f$ with respect to $u$.
We shall need this expansion up to the third order. Denoting derivatives with respect to u by primes, we have
\begin{equation}  
\phi(z, u+iP) = \phi(z,u) + i \phi'(z,u)P - \frac12
\phi''(z,u)P^2 - \frac{i}6\phi'''(z,u)P^3 + \dots\ \tag{3.12}
\end{equation}
In (3.9) we denote
$$L(f,g) =  Re \{ig(z,u+i\pz) +
2 P_z f(z,u + i\pz)\}.$$
From (3.12) we have
\begin{equation}
\begin{aligned}
  2L(f,g)\  =&2 f P_z +2i f' P_zP  - f'' P_zP^2 +
2 \bar f \bar P_z -2i \bar f' \bar P_z \bar P  - \bar f'' \bar P_z\bar P^2 \\
  +   i g & - g' P - \frac{i}2 g'' P^2 + \frac{1}6 g''' P^3 -
 i \bar g - \bar g' \bar P + \frac{i}2 \bar g''
\bar P^2 + \frac{1}6 \bar g''' \bar P^3  + \dots\
\end{aligned}
 \tag{3.13} 
\end{equation}
In this expansion we will collect terms of type
 $(i,j).$

In order to compute  the coefficients of
 terms of the types
specified by (3.11), let us  denote  $F^*_{ij} = a^*_{ijm}$ and
$F_{ij} = a_{ijm}$ if
 $i + j + km = \mu$,
$f_i = f_{ij}$ if $ i + kj = \mu - k + 1$ and $g_i = g_{ij} $ if
$i + kj = \mu$. Similarly, for the derivatives of $f(z,u)$ and
$g(z,u)$ with respect to $u$, we denote $f'_i= j f_{ij}$, $f''_i =
j(j-1)f_{ij}$ if  $ i + kj = \mu - k + 1$ and $g'_i= j g_{ij}$,
$g''_i = j(j-1)g_{ij}$, $g'''_i = j(j-1)(j-2)g_{ij}$,
  if  $ i + kj = \mu$.
To separate this notation from the one used for weights,
subscripts indicating weight are always greek letters.
 Finally,
in analogy to (3.10) we write

\begin{equation}
(F_{k-1}, P_z) = \sum_{j=1}^{k-2}F_{j,k-1-j} (j+1)\bar a_{j+1}, \ \tag{3.14}
\end{equation}
and the same for $F^*$. For $j\ge 1$ and $j\neq k-1$ we have from
(3.9) and (3.13)
\begin{equation}
F^*_{j0} = -\frac i2 \ g_j + F_{j0} + \dots \ \tag{3.15}
\end{equation}
 Here
dots denote terms with values already determined, i.e. terms
depending on $f_{\nu-k+1}$, $g_{\nu}$, $F_{\nu}$, $F^*_{\nu}$ for $\nu <
\mu$.
 For
$j=k-1$ we have
\begin{equation}
F^*_{k-1,0} = -\frac i2 \ g_{k-1}   - \bar f_0
\delta_{1l}+ F_{k-1,0} + \dots ,\ \tag{3.16} 
\end{equation}
 where
$\delta_{1l}$ is Kronecker's delta. Further, for $j \geq 1$ and $j
\neq k-1$ only the first and the eighth term in (3.13) contribute
to type $(k-l+j,l)$, and we get
\begin{equation}
 F^*_{k-l+j,l}
 = -(k-l) f_{j+1} + \frac12  g'_j + F_{k-l+j, l} + \dots \ \tag{3.17}
\end{equation}
For $j=k-1$ we have
\begin{equation}
F^*_{2k-l-1,l} = -(k-l)f_{k} + \frac12 g'_{k-1} + i \bar f_0
\delta_{1l} + \dots \ \tag{3.18}
\end{equation}
 Next
\begin{equation}
F^*_{k-l, l}  = -(k-l) f_{1} - l \bar f_1 +  Re g'_0 + F_{k-l,l} +
\dots \ \tag{3.19}
\end{equation}
 and
\begin{equation}
F^*_{2k-2l,2l} = -(k-2l) f_{k+1} \bar a_{2l} - i (k-l) f'_1 + i l \bar f'_1
 + \frac12 g'_k \bar a_{2l} - \frac12 Im g''_0
+ F_{2k-2l,2l} + \dots \ \tag{3.20}
\end{equation}
 Further,
\begin{equation}
(F^*_{k-1}, P_z) = f_0 (P_z,P_z) + \bar f_0 (P_{\bar z}, P_z) +
 (F_{k-1}, P_z) + \dots \ \tag{3.21}
\end{equation}
We first show that (3.21) determines uniquely $f_0$. For this we
use the elementary fact that an equation $\alpha f_0 + \beta \bar
f_0 = \gamma$, where $\alpha, \beta, \gamma \in \mathbb C$,
determines uniquely $f_0$ if and only if $\vert \alpha \vert \neq
\vert \beta \vert$.
 Since  $\bar{P_{\bar z}} = P_z$, we have  $   (P_z, P_z) =( P_{\bar
 z}, P_{\bar z})$. By Cauchy-Schwartz inequality
$$\vert  (P_z, P_{\bar z}) \vert \leq \lbrace  (P_{\bar z},P_{\bar z} ) (P_z, P_z)
 \rbrace^{\frac12} =  (P_z, P_z),$$
 with equality
if and only if $P_z$ is a scalar multiple of $P_{\bar z}$, i.e., $
P_z = e^{i\theta} P_{\bar z}$, modulo harmonic terms. It is easy
to verify that this happens if and only if $P$  is equivalent to
 $(Re\;z)^k$, the exceptional case treated in section 4.
Next, (3.15) and (3.16) determine $g_j$ for all $j \geq 1$. Then
(3.17) and (3.18) determine $f_j$ for all $j\geq 2$. The pair of
equations (3.19), (3.20) then determines $f_1$ and $g_0$. It is
straightforward to verify that the initial appearence of each
equation agrees with our normalization conditions for $f$ and $g$.

\section{Normal forms for exceptional models}

Now we consider the remaining two cases.  In the first one, when
$2l = k$,
 the model hypersurface is described by
$$ v = \ab k.$$
It is preserved by transformations of the form
\begin{equation}
z^* =  \frac{ \delta e^{i\theta} z}{(1 + \mu w)^{\frac1l}}, \ \  \ \ \
  w^* =\frac{ \delta^k w}{1 + \mu w},\tag{4.1}
\end{equation}
                           where $\delta > 0,$ and $\te, \mu  \in \mathbb R$.
We denote this group again by $H$. Its real dimension is equal to
three. We consider hypersurfaces of the form
$$ v = \ab k + \text{terms of weight} \ge  k+1.$$
Now we define  $\mathcal F$ as in section 3. Again, $\mathcal T_0$
will denote the group of transformations which preserve this form.
Its elements are of the form
\begin{equation}
z^* =\delta e^{i\theta} z + \text{terms of weight} \ge 2,\ \  \ \ \  
 w^* = \delta^k w + \text{terms of weight}
 \ge k+1,\tag{4.2}
\end{equation}
where $\theta \in \mathbb R$ and $\delta \in \mathbb R^+$.

 $\mathcal T$ will denote the group of transformations of the form (3.7),
 satisfying
an additional  normalization condition
$$Re \; g_{ww} = 0.$$
It is verified easily that we can factor any  $\tau \in \mathcal
T_0$ uniquely as $\tau = \phi \circ T$, where $ \phi \in H $ and
$T\in \mathcal T$. Indeed, the additional  condition determines
$\mu$ in (4.1), while $\delta$ and $e^{i\theta}$ are determined by
(3.8).

{\bf Definition 4.1.} We say that  $F$ is in  normal form if
$$
\begin{aligned}
Z_{j0}&= 0, \ \ \ \ \ j=0,1,\dots,  \\
Z_{l,l+j}&= 0, \ \ \ \ \ j= 1,2,\dots, \\
Z_{l,l}& = 0, \\
Z_{2l, 2l}& = 0,\\
Z_{3l, 3l}& = 0,   \\
Z_{2l, 2l-1}& = 0. \ \end{aligned} 
$$

 We have

{\bf Proposition 4.2.} {\it For any $F \in \mathcal F$ there is a
uniquely determined formal transformation $T \in \mathcal T$ which
transforms $F$ into normal form.}

{\it{proof. }} With the notation introduced in the proof of
Proposition 3.3, we obtain for $j \geq 0$ from (3.9) and (3.13)
\begin{equation}
 \begin{aligned}
F^*_{j0}& = -\frac i2 \ g_j + F_{j0} + \dots \\
F^*_{00}& =  Im \ g_0 + F_{00} + \dots \ 
\end{aligned}
\tag{4.3}
\end{equation}
For $j \geq 1$ we have
\begin{equation}
\begin{aligned}
 F^*_{ll}&=  Re\ g'_0 - 2l Re \ f_1 + F_{ll} + \dots  \\
F^*_{l+j,l}& = - l f_{j+1} + \frac12  g'_j + F_{l+j,l} + \dots \\
F^*_{2l,2l}& = -\frac12 Im \ g''_0  + 2l Im\ f'_1 + F_{2l,2l} + \dots \\
F^*_{3l,3l}& = -\frac16 Re \ g'''_0 + lRe\ f''_1 + F_{3l,3l} + \dots \\
F^*_{2l, 2l-1}& =-  il \bar f'_0 + F_{2l, 2l-1} + \dots \
\end{aligned}
\tag{4.4} 
\end{equation}
Again, the  condition that $F^*$  be in normal form determines
uniquely coefficients of $f$ and $g$. Equations (4.3) determine
$g_j$ for all $j\geq 1$ and $Im \ g_0$. The second  equation in
(4.4) determines $f_j$ for $j \geq 2$ and the fifth $f_0$. The
third equation determines $Im \ f_1 $.
  Finally the first and the fourth equations  determine $Re \ f_1$ and
$Re \ g_0$.

Now we consider the model hypersurface equivalent to $\{v = (Re \;
z)^k\}$. In this case
$$P(z,\bar z) = \frac1k [(z + \bar z )^k - 2Re\; z^k].$$

{\bf Definition 4.3.}
 We say that   $F$ is in  normal form if
\begin{equation}
\begin{aligned} Z_{j0} &= 0, \ \ \ \ \ j=1,2,\dots,  \\
Z_{k-l+j,l} &= 0, \ \ \ \ \ j= 0,1,\dots, \\
Z_{2k-2l, 2l} & = 0, \\
Re \; Z_{k-2,1}  & = 0, \\
Re \; Z_{k, k-1} & = 0.
 \ 
\end{aligned}
\tag{4.5}
\end{equation}

$\mathcal T_0$ and $\mathcal T$ are now defined as for generic
models in part 3. Now we prove Proposition 4.2 for this case.

{\it proof of Proposition 4.2} We proceed as in the proof of
Proposition 3.3. Instead of (3.21) we obtain
$$ Re \; F^*_{k-2,1} = -2(k-1) Re \;f_0 + Re \; F_{k-2,1} + \dots$$
and
$$Re \; F^*_{k, k-1} = (2C-1) Im\; f_0  + \frac12 Re\;g'_{k-1} + Re \; F_{k,k-1}, $$
where $C = \frac1k [\binom{2k-1}k-1] $ is the coefficient of $z^k \bar z^{k-1}$
 in $P P_{\bar z}$.
Since $k \ge 3 $, we have $C \geq 3$, and  this pair of equations together with
  (3.16) determine uniquely $g_{k-1}$ and $f_0$.

\section{The  equivalence problem and applications}

Propositions 3.3 and 4.2 combined with the result of [BER] give a
solution to the local equivalence problem. Let $\mathcal N$ be the
subset of $\mathcal F$ containing  formal power series in normal
form. Let $M$ be a hypersurface described by (3.6) (where the
exceptional cases of P are now included),
 and let $T_0 \in \mathcal T_0$
be a  transformation which takes $M$ into normal form. Writing
$T_0 = \phi \circ T$, with $\phi \in H$ and $ T \in \mathcal T$,
we see from Propositions 3.3 and 4.2 that transformations into
normal form are
 parametrized by elements of $H$. If $M$ is already in normal form,
we get a natural action of $H$ on $\mathcal N$.
 It is obtained by
first applying an element of $H$ and then renormalizing by a
unique element in $\mathcal T$.

  Since normal forms are unique up to
this action of $H$, we consider the set of equivalence classes
$\mathcal N \ \text{mod}\ H$ and call the elements of this set
normal forms.

 {\textbf{Theorem 5.1.}}
{\it Two real analytic hypersurfaces are
locally biholomorphically equivalent if and only if their normal
forms are equal.
}

 {\it{proof. }}
Clearly the normal form is a biholomorphic invariant. On the other
hand, if $\phi_1$, $\phi_2$ are formal transformations of $M_1$,
$M_2$ respectively into the same power series in normal form, then
$\phi_1 \circ \phi_2^{-1}$ is a formal equivalence of $M_1$,
$M_2$. Since the hypersurfaces are of finite type, the result of
[BER] implies that it converges. \

For nondegenerate hypersurfaces it may be quite difficult to
decide \newline whether two  hypersurfaces in normal form are equivalent
under the action of $H$, which has real dimension five. In the
degenerate case, this becomes much simpler. As an example, let us
consider the following case.
   Let $M_1, M_2$ be two hypersurfaces of finite type $k$ in normal form,
which for
$M_1$ is given by
$$ P(z, \bar z) + \sum_{i+j < k-1} \sum_m a_{ijm}z^i \bar z^j u^m$$
and for $M_2$ by
 $$ P(z, \bar z) + \sum_{i+j < k-1} \sum_m b_{ijm}z^i \bar z^j u^m,$$
where $P$ is a polynomial of the form (2.10), (3.4), (3.5) and the
sums on the right contain only terms of weight bigger than $k$. We
have the following consequence of Proposition 3.3.

{\textbf{Proposition 5.2.} }
{\it $ M_1$ and $M_2$ are  locally
biholomorphically equivalent if and only if there is  $(\delta ,
e^{\sqrt{-1}\theta})\in H$ such that
$$a_{ijm} = \delta^{i+j+km} e^{\sqrt{-1}(i-j)\theta } b_{ijm}$$ for
 all indices $(i,j,m)$.
}

As a corollary of Propositions 3.3 and 4.2 we obtain

{\textbf{Corollary 5.3.}}
{\it    The only transformations which
preserve the model hypersurface are the elements of $H$.
}

{\it{proof. }}  Let $\tau \in \mathcal T_0$ be a transformation
which preserves the
 model hypersurface. There is a unique $\phi \in H$ such that
$\phi \circ \tau$ is an element of $\mathcal T$. This mapping
still preserves the model hypersurface, in particular it is a
mapping into normal form. By Propositions 3.3 and 4.2, $\phi \circ
\tau$ is the identity, therefore $\tau$ is an element of $ H$

Propositions 3.3 and 4.2 also give immediately precise information
about the dimension
 of the  stability group (the local symmetry group)
of $M$ at $p$.

{\textbf{Corollary 5.4.}}
{\it
 Let $M$ be a real analytic hypersurface
in $\mathbb C^2$ and let $p \in M$ be a point of finite type $k$,
where $k \geq 3$. Then the dimension of the stability group of $M$
at $p$ is less or equal to three.
If, moreover, the model hypersurface at $p$ is different from $\{
v = \vert z \vert^k\}$, then the dimension is at most one.
}

\vspace{5mm}
\begin{center}
{\textsc{Acknowledgement:}} \end{center}
 {Supported by 
 a grant of the GA \v CR no. 201/05/2117}
\end{document}